\documentclass[11pt]{article}
\usepackage{amsmath,latexsym,amsbsy,amssymb,amsthm,color,enumitem}
\usepackage{graphicx}
\usepackage{graphicx}

\def\ds{\displaystyle}
\def\forall{\hbox{for all}~}

\def\ve{\varepsilon}

\def\R{\mathbb{R}}

\def\vs{\vskip 2em}

\def\v{\vskip 1em}

\def\bega{\begin{array}}
\def\enda{\end{array}}
\def\begi{\begin{itemize}}
\def\endi{\end{itemize}}

\def\bel{\begin{equation}\label}
\def\eeq{\end{equation}}
\def\sqr#1#2{\vbox{\hrule height .#2pt
\hbox{\vrule width .#2pt height #1pt \kern #1pt
\vrule width .#2pt}\hrule height .#2pt }}
\def\square{\sqr74}
\def\endproof{\hphantom{MM}\hfill\llap{$\square$}\goodbreak}

\newtheorem{theorem}{Theorem}[section]
\newtheorem{corollary}[theorem]{Corollary}
\newtheorem{definition}[theorem]{Definition}
\newtheorem{remark}[theorem]{Remark}
\newtheorem{lemma}[theorem]{Lemma}
\newtheorem{proposition}{Proposition}[section]

\begin{document}
\title{\bf Metric entropy for  functions of bounded total generalized variation}\vs
\author{\it Rossana Capuani, Prerona Dutta and Khai T. Nguyen}
\author{Rossana Capuani, Prerona Dutta, and Khai T. Nguyen\\ 
\\
  {\small  Department of Mathematics, North Carolina State University }\\ 
 \\  {\small e-mails: ~rcapuan@ncsu.edu, ~pdutta2@ncsu.edu,  ~ khai@math.ncsu.edu}
 }

\maketitle
\begin{abstract}
We  establish a sharp estimate for a minimal number of binary digits (bits) needed to represent all  bounded total generalized  variation functions taking values in a general totally bounded metric space $(E,\rho)$ up to an accuracy of $\ve>0$ with respect to the ${\bf L}^1$--distance. Such an estimate is explicitly computed in terms of  doubling and packing dimensions of $(E,\rho)$. The obtained result is applied to provide an upper bound on the metric entropy for a set of entropy admissible weak  solutions to scalar conservation laws in one-dimensional space with  weakly genuinely nonlinear fluxes.
\end{abstract}
\section{Introduction}
\label{sec:1}
\setcounter{equation}{0}
The metric entropy (or ${\ve}$-entropy) has been studied extensively in a variety of literature and disciplines.  It plays a central role in various areas of information theory and statistics, including nonparametric function estimation, density information, empirical processes and machine learning (see e.g in \cite{LB, DH, DP}). It provides a tool for characterizing the rate of mixing of sets of small measure. The notion of metric entropy (or $\ve$-entropy) has been introduced by Kolmogorov and Tikhomirov \cite{KT} in 1959 as follows:
\medskip

\begin{definition}\label{entropy} 
Let $(E,\rho)$ be a metric space and $K$ be a totally bounded subset of $E$. For $\varepsilon >0$, let $\mathcal{N}_{\varepsilon}(K\big| E)$ be the minimal number of sets in an $\varepsilon$-covering of $K$, i.e.,  a covering of $K$ by balls in $E$ with radius no greater than $\varepsilon$.
Then the $\varepsilon$-entropy of $K$ is defined as 
\[
\mathcal{H}_{\varepsilon}(K|E)~=~\log_2\mathcal{N}_{\varepsilon}(K|E).
\]
\end{definition}
 A classical topic in the field of probability is to investigate the metric covering numbers for general classes  $\mathcal{F}$ of real-valued functions defined on $E$ under the family of ${\bf L}^1(dP)$ where $P$ is a probability distribution on $E$. Upper and lower bounds on the ${\ve}$-entropy of $\mathcal{F}$ in terms of Vapnik-Chervonenkis, pseudo-dimension and the scale-sensitive dimension of the function class were established in \cite{RMD,DH, DH1, KMJ, DP} and in \cite{KMJ, LBW}.
\medskip

Thanks to the Helly's theorem, a set of uniformly bounded variation functions is compact in ${\bf L}^1$-space.  Consequently, attempts were made  to quantify the degree of compactness of such sets by using the $\ve$-entropy. In \cite{KMJ}, the authors showed that  the $\ve$-entropy of any set of uniformly bounded total variation real-valued functions in ${\bf L}^1$ is of the order $\ds {1\over \ve}$ in the scalar case. Later on, this result was also extended to  multi-dimensional cases in \cite{DN}. Some related works have been done in the context of density estimation where attention has been given to the problem of finding covering numbers for the classes of densities that are unimodal or non-decreasing  in \cite{LB, PG}. In the multi-dimensional cases, the covering numbers of convex and uniformly bounded functions  were studied in \cite{GS}. It was shown that the $\ve$-entropy of a class of convex functions with uniform bound in ${\bf L}^1$ is of the order $\ds{1\over \ve^{d\over 2}}$ where $d$ is the dimension of the state variable. The result was previously studied for scalar state variables in \cite{DD} and for convex functions that are uniformly bounded and uniformly Lipschitz with a known Lipschitz constant in \cite{EB}. These results have direct implications in the study of rates of convergence of empirical minimization procedures (see in \cite{LB1,SVG}) as well as optimal convergence rates in the numerous convexity constrained function estimation problems (see in \cite{LB0, LLC, YB}).
\medskip

From a different aspect, the ${\ve}$-entropy has been used to measure the set of solutions of  nonlinear partial differential equations. In this setting, it  could provide a measure of the order of ``resolution'' and the ``complexity'' of a numerical scheme, as suggested in~\cite{Lax78}. The first results on this topic  were obtained in \cite{AON1, DLG} for  the scalar conservation law with uniformly convex flux $f$ (i.e. $f''(u)\geq c>0$),  in one-dimensional space
\bel{1}
u_t(t,x)+f(u(t,x))_x~=~0\,.
\eeq
It was shown  that the number of functions needed to represent an entropy admissible weak solution $u$ at any time $t>0$ with an accuracy of $\varepsilon$ with respect to the $\bf{L}^1$-distance is of the order  $\ds {1\over\ve }$. A similar estimate was also obtained  for the system of hyperbolic conservation laws in \cite{AON2, AON3} and for Hamilton-Jacobi equations with uniformly convex Hamiltonian in \cite{ACN, ACN1}.  All these proofs strongly relied on the BV regularity properties of solutions. Thereafter, the results in \cite{AON1, DLG} were extended to scalar conservation laws with a smooth flux function $f$ that is either strictly (but not necessarily uniformly) convex or has a single inflection point with a polynomial degeneracy \cite {AON4} where entropy admissible weak solutions may have unbounded total variation.  In this case, the sharp estimate on the $\ve$-entropy for sets of entropy admissible weak solutions was provided by exploiting the BV bound of  the characteristic speed $f'(u)$ at any positive time \cite{KSC86}. On the other hand, it was shown  in \cite[Example 7.2]{SM}) that for fluxes having one inflection point where all derivatives vanish, the composition of the derivative of the flux with the solution of (\ref{1}) fails in general to belong to the BV space and the analysis in \cite{AON4} cannot be applied here. However, for {\it weakly genuinely nonlinear} fluxes, that is to say for fluxes with no affine parts,  equibounded sets of entropy solutions  of (\ref{1}) at positive time are still relatively compact in ${\bf L}^1$ (see \cite[Theorem 26]{Ta}). Therefore, for fluxes of such classes that do not fulfill the assumptions in \cite{AON4}, it remains an open problem to provide a sharp estimate on the $\ve$-entropy for the  solution set of (\ref{1}). A different approach from \cite{AON4} must be pursued to study the $\ve$-entropy for (\ref{1}) with weakly genuinely nonlinear fluxes, perhaps exploiting the uniform bound on total generalized variation  of entropy admissible weak solutions studied in~\cite[Theorem 1]{Marc}.
\medskip

From the above viewpoints, the present paper aims to study the $\ve$-entropy of classes of uniformly bounded total generalized  variation functions taking values in a general totally bounded metric space $(E,\rho)$. More precisely, for a given convex function $\Psi: [0,+\infty)\to [0,+\infty)$ with $\Psi(0)=0$ and $\Psi(s)>0$ for all $s>0$, let $\mathcal{F}^{\Psi}_{[L,V]}$ be a set of functions $g: [0,L]\to E$ such that the $\Psi$-total variation of $g$ over the interval $[0,L]$ is bounded by $V$, i.e., 
\[
\sup_{N\in\mathbb{N}, 0=x_0<x_1<...<x_N=L}\sum_{i=0}^{N-1}\Psi \left(\rho\left(g(x_{i}), g(x_{i+1}\right))\right)~\leq~V.
\]
We establish  upper and  lower bounds on $\mathcal{H}_{\ve}\left(\mathcal{F}^{\Psi}_{[L,V]}\Big|{\bf L}^1([0,L],E)\right)$, the $\ve$-entropy of $\mathcal{F}^{\Psi}_{[L,V]}$ with respect to the ${\bf L}^1$-distance. For deriving sharp estimates explicitly, our idea is to  use the notions of doubling and packing dimensions of $(E,\rho)$, {\color{black} denoted by ${\bf d}(E)$ and ${\bf p}(E)$ respectively, which were first introduced by Assouad in \cite{PA}.  In Theorem \ref{main}, we prove that for every $\ve>0$ sufficiently small, the sharp bounds on $\mathcal{H}_{\ve}\left(\mathcal{F}^{\Psi}_{[L,V]}\Big|{\bf L}^1([0,L],E)\right)$ can be approximated in terms of ${\bf p}(E)$, ${\bf d}(E)$ and $\Psi$. In particular, if  $\Psi(s)=s^{\gamma}$ for some $\gamma\geq 1$ and the metric space $(E,\rho)$ is generated by a finite dimensional normed space $(\R^d,\|\cdot\|)$ then  the $\ve$-entropy of $\mathcal{F}^{\Psi}_{[L,V]}$ in ${\bf L}^1\left([0,L],\R^d\right)$ is of the order $\ds {d\over \ve^{\gamma}}$, i.e., 
\[
\mathcal{H}_{\ve}\left(\mathcal{F}^{\Psi}_{[L,V]}\Big|{\bf L}^1\left([0,L],\R^d\right)\right)~\approx~\ds {d\over \ve^{\gamma}}~.
\]
The result is applied to provide an upper estimate on the $\ve$-entropy  of a set of entropy admissible weak solutions to scalar conservation laws (\ref{1}) with general weakly genuinely nonlinear fluxes in Theorem \ref{M1}, which partially extends the recent one in \cite{AON4}.  The estimate is sharp in the case of fluxes having finite  inflection points with a polynomial degeneracy.  However, a natural question regarding sharp estimates of the $\ve$-entropy for such solution sets to (\ref{1}) with general weakly genuinely nonlinear fluxes is still open.
\medskip 

This paper is organised as follows. In Section 2, we present some preliminary results on covering and packing numbers of a totally bounded metric space and also include necessary concepts related to functions of bounded total generalized variation. In Section 3, the first subsection focuses on finding the upper and lower estimates of the $\ve$-entropy for a set of bounded total generalized variation functions, while the second subsection is an application of these estimates to scalar conservation laws with weakly genuinely nonlinear fluxes.
 \section{Notations and preliminaries}
\setcounter{equation}{0}
Let $E$ be a metric space with distance $\rho$ and $I$ be an interval in $\R$. Throughout the paper we shall denote by:
\begin{itemize}
	\item $B_\rho(z,r)$, the open ball of radius $r$ and center $z$, with respect to the metric $\rho$ on $E$, i.e.,
	\[
	B_\rho(z,r)~=~\left\{y\in E~|~\rho(z,y)<r\right\};
	\]
	\item $\mathrm{diam}(F)=\sup_{x,y\in F} \rho(x,y)$, the diameter of the set $F$ in $(E,\rho)$;
	
	\item $\mathbf{L}^{1}(I,E)$, the Lebesgue metric space of all (equivalence classes of) summable functions $f: I\to E$, equipped with the usual ${\bf L}^1$-metric distance, i.e., 
	\[
	\rho_{{\bf L}^1}(f,g)~:=~~\int_{I}\rho(f(t),g(t))dt~<~+\infty
	\]
	for every $f,g\in \mathbf{L}^{1}(I,E)$;
	\item ${\bf L}^1(\mathbb{R})$, the Lebesgue space of all (equivalence classes of) summable functions on $\mathbb{R}$, equipped with the usual norm $\| \cdot\|_{{\bf L}^1}$;
	\item ${\bf L}^\infty(\mathbb{R})$, the space of all essentially bounded functions on $\mathbb{R}$, equipped with the usual norm $\| \cdot\|_{{\bf L}^\infty}$;
	\item $\mathrm{Supp} (u)$, the essential support of a function $u\in {\bf L}^\infty(\mathbb{R})$;
	
	\item $B_{{\bf L}^1(I,E)}(\varphi, r)$, the open ball of radius $r$ and center $\varphi$ in $\mathbf{L}^{1}(I,E)$,  with respect to the metric $\rho_{{\bf L}^1}$ on $\mathbf{L}^{1}(I,E)$, i.e.,
	\[
	B_{{\bf L}^1(I,E)}(\varphi, r)~=~\left\{g\in \mathbf{L}^{1}(I,E)~\big|~\rho_{{\bf L}^1}(\varphi,g)<r\right\};
	\]
	
	\item $\mathcal{B}(I,[0,+\infty))$, a set of bounded functions from $I$ to $[0,+\infty)$;
	\item $\mathcal{C}^{\infty}(\R,\R)$, space of smooth functions having derivatives of all orders;
	\item $TV(g,I)$, total variation of $g$ over the interval $I$;
	\item $TV^{\Psi}(g,I)$, $\Psi$-total variation of $g$ over the interval $I$;
	\item $TV^{{1\over \gamma}}(g,I)$, $\gamma$-total variation of $g$ over the interval $I$, i.e., $\Psi$-total variation of $g$ with $\Psi$ defined by $\Psi (s) = |s|^{\gamma}$;
	\item $\chi_I(x)=\left\{\bega{rl}
	&1 \qquad~~\mathrm{if}\qquad x\in I\,
	\\[4mm]
	&0\quad~\mathrm{if}\qquad x\in\R^n\backslash I\,
	\enda\right.$
	the characteristic function of $I$;
	\item $\mathrm{Card}(S)$, the number of elements in any finite set $S$;
	\item $\lfloor x\rfloor := \max\{z\in\mathbb{Z}~|~z\leq x\}$, the integer part of $x$;
	\item $\overline{1,N}$, the set of natural numbers from $1$ to $N$;
	\item $\ds{n \choose k} = {n!\over k! (n-k)!}$, number of ways in which $k$ objects can be chosen from among $n$ objects.
\end{itemize}
\subsection{Covering, packing and metric dimension}
Let us first recall the concepts of covering number and packing number in a totally bounded metric space $(E,\rho)$. For any $K\subseteq E$ and $\alpha>0$, we say that 
\begin{itemize}
\item the set $\mathcal{A} = \{a_1,a_2,\dots, a_n\}\subseteq E$ is an $\alpha$-covering of $K$ if $K\subseteq \bigcup_{i=1}^{n}B_\rho(a_i,\alpha)$, or equivalently, for every $x\in K$, there exists $i\in \overline{1,n}$ such that $\rho(x,a_i)<\alpha$; $\mathrm{Card}(\mathcal{A})$ is called the {\it size} of this $\alpha$-covering;
\item the set $\mathcal{B} = \{b_1,b_2,\dots, b_m\}\subseteq K$ is an $\alpha$-packing of $K$ if $\rho(b_i,b_j)>\alpha$ for all $i\neq j\in \overline{1,m}$, or equivalently, $\{B_\rho(b_i,\alpha/2)\}^{m}_{i=1}$ is a finite set of disjoint balls; $\mathrm{Card}(\mathcal{B})$ is called the {\it size} of this $\alpha$-packing.
\end{itemize}
\v
\begin{definition}\label{pc} 
The $\alpha$-covering and $\alpha$-packing numbers of $K$ in $(E,\rho)$ are defined by 
\[
\mathcal{N}_{\alpha}(K|E)~=~\min\left\{n\in\mathbb{N}~|~\exists~\alpha\mathrm{-covering~of}~K~\mathrm{having~size}~n\right\}
\]
and
\[
\mathcal{M}_{\alpha}(K|E)~=~\max\left\{m\in\mathbb{N}~|~\exists~\alpha\mathrm{-packing~of}~K~\mathrm{having~size}~m\right\},
\]
respectively.
\end{definition}
Since  $E$ is totally bounded, $\mathcal{N}_{\alpha}(K|E)$ is finite for every $\alpha>0$. Moreover, the maps $\alpha\mapsto \mathcal{N}_{\alpha}(K|E)$ and $\alpha\mapsto \mathcal{M}_{\alpha}(K|E)$ are non-increasing. The relation between $\mathcal{N}_{\alpha}(K|E)$ and $\mathcal{M}_{\alpha}(K|E)$ is described by the following simple double inequality: 
\medskip

\begin{lemma}\label{Rl} For any $\alpha>0$, one has
\[
\mathcal{M}_{2\alpha}(K|E)~\leq~\mathcal{N}_{\alpha}(K|E)~\leq~\mathcal{M}_{\alpha}(K|E).
\]
\end{lemma}
{\bf Proof.}
For the proof see e.g in \cite{KT}.
\endproof
\medskip

Let us now introduce a commonly used notion of dimension for a metric space $(E,\rho)$, as proposed in \cite[\S 4]{PA}.
\v
\begin{definition}\label{pc-dimension} The doubling and  packing dimensions of $(E,\rho)$ are respectively defined by
\begin{itemize}
\item ${\bf d}(E)$ is the minimum natural number $n$ such that  for every $x\in E$ and $\alpha>0$, the ball $B_\rho(x,2\alpha)$ can be covered by $2^{n}$ balls of radius $\alpha$;
\item ${\bf p}(E)$  is the maximum natural number $m$ such that  for every $x\in E$ and $\alpha>0$, the ball $B_\rho(x,2\alpha)$ contains an $\alpha$-packing of size $\mathcal{M}_{\alpha}(B_\rho(x, 2\alpha)|E)$ which satisfies the inequality 
$$2^{m} ~\leq~ \mathcal{M}_{\alpha}(B_\rho(x, 2\alpha)|E)~<~2^{m+1}.$$
\end{itemize}

\end{definition}
We conclude this subsection with the following result relating $\alpha$-covering and $\alpha$-packing.
\medskip
\begin{lemma}\label{co-R-d}
Given $R\geq 2\alpha >0$, let $k$ and $m$ be   natural numbers  such that
\[
2\cdot 7^{k}~\leq~ {R\over \alpha}~\leq~2^{m}.
\]
The following hold
\bel{ese1}
\mathcal{N}_{\alpha}\left(B_\rho(z,R)~\Big|~E\right)~\leq~2^{m{{\bf d}(E)}}
\eeq
and 
\bel{esee1}
\mathcal{M}_{\alpha}\left(B_\rho(z,R)~\Big|~E\right)~\geq~2^{(k+1) {\bf p}(E)}
\eeq
for all $z\in E$.
\end{lemma}
{\bf Proof.} {\bf 1.} For every $n\geq 0$, we first show that 
\bel{Ne1}
\mathcal{N_{\alpha}}\left(B_\rho(z,2^{n}\alpha)~\Big|~E\right)~\leq~2^{n{{\bf d}(E)}}\qquad\forall  z\in E.
\eeq
Assume that (\ref{Ne1}) holds for $n=i\geq 0$. For any given  $z_0\in E$, from Definition \ref{pc-dimension}, one has 
\[
\mathcal{N}_{2^i\alpha}\left(B_\rho(z_0,2^{i+1}\alpha)~\Big|~E\right)~\leq~2^{{\bf d}(E)}.
\]
Equivalently, there exist $x_1,x_2,\dots, x_{2^{{\bf d}(E)}}\in E$ such that 
\[
B_\rho(z_0,2^{i+1}\alpha)~\subseteq~\bigcup^{2^{{\bf d}(E)}}_{j=1}B_\rho(x_j,2^{i}\alpha)
\]
and   
\[
\mathcal{N_{\alpha}}\left(B_\rho(z_0,2^{i+1}\alpha)~\Big|~E\right)~\leq~\sum_{j=1}^{2^{{\bf d}(E)}}\mathcal{N_{\alpha}}\left(B_\rho(x_j,2^{i}\alpha)~\Big|~E\right)~\leq~2^{{\bf d}(E)}\cdot 2^{i{\bf d}(E)}~=~2^{(i+1){\bf d}(E)}.
\]
Thus, (\ref{Ne1}) holds for $n=i+1$ and the method of induction yields (\ref{Ne1}) for all $n\geq0$. In particular, the non-decreasing property of  the map  $r\mapsto\mathcal{N_{\alpha}}\left(B_\rho(z,r)~\Big|~E\right)$ implies that 
\[
\mathcal{N}_{\alpha}\left(B_\rho(z,R)~\Big|~E\right)~\leq~\mathcal{N_{\alpha}}\left(B_\rho(z,2^{m}\alpha)~\Big|~E\right)~\leq~2^{m{{\bf d}(E)}}.
\]

{\bf 2.} To achieve the  inequality in (\ref{esee1}), we prove that 
\bel{Me1}
\mathcal{M}_{\alpha}\left(B_\rho(z,2\cdot 7^{n}\alpha)~\Big|~E\right)~\geq~2^{(n+1){\bf p}(E)}\qquad\forall z\in E.
\eeq
It is clear from Definition \ref{pc-dimension} that (\ref{Me1}) holds for $n=0$. Assume that (\ref{Me1}) holds for $n=i\geq 1$. For any given  $z_0\in E$, from Definition \ref{pc-dimension}, one has 
\[
\mathcal{M}_{6\cdot 7^{i}\alpha}\left(B_\rho(z_0,12\cdot 7^{i}\alpha)~\Big|~E\right)~\geq~2^{{\bf p}(E)}.
\]
Equivalently, there exist $x_1,x_2,\dots, x_{2^{{\bf p}(E)}}\in B_\rho(z_0,12\cdot 7^{i}\alpha)$ such that 
\[
\rho(x_{j_1},x_{j_2})~>~6\cdot 7^{i}\alpha~\geq~4\cdot 7^{i}\alpha+2\alpha \qquad\forall j_1\neq j_2\in \{1,2,\dots, 2^{{\bf p}(E)}\}.
\]
In particular, for every  $j_1\neq j_2\in \{1,2,\dots, 2^{{\bf p}(E)}\}$, it holds
\[
\rho(z_1,z_2)~>~2\alpha\qquad\forall z_1\in B_\rho(x_{j_1},2\cdot 7^{i}\alpha), z_2\in B_\rho(x_{j_2},2\cdot 7^{i}\alpha).
\]

Since $B_\rho(x_{j},2\cdot 7^{i}\alpha)\subseteq B_\rho(z_0,2\cdot 7^{i+1}\alpha)$ for all $j\in  \{1,2,\dots, 2^{{\bf p}(E)}\}$, one then has
\begin{eqnarray*}
\mathcal{M}_{\alpha}\left(B_\rho(z_0,2\cdot 7^{i+1}\alpha)~\Big|~E\right)&\geq&\sum_{j=1}^{2^{{\bf p}(E)}}\mathcal{M}_{\alpha}\left(B_\rho(x_j,2\cdot 7^{i}\alpha)~\Big|~E\right)~\geq~2^{{\bf p}(E)}\cdot 2^{(i+1){\bf p}(E)}\\
&=&2^{(i+2){\bf p}(E)}.
\end{eqnarray*}
Thus, by the method of induction, (\ref{Me1}) holds for all $n\geq 0$. In particular, the non-decreasing property of  the map  $r\mapsto\mathcal{M_{\alpha}}\left(B_\rho(z,r)~\Big|~E\right)$ implies that 
\[
\mathcal{M}_{\alpha}\left(B_\rho(z,R)~\Big|~E\right)~\geq~\mathcal{M_{\alpha}}\left(B_\rho(z,2\cdot7^k\alpha)~\Big|~E\right)~\geq~2^{(k+1){{\bf p}(E)}}.
\]
\qed
\medskip

As a consequence of Lemma \ref{Rl} and Lemma \ref{co-R-d}, one has that 
\bel{ese3}
\left({R\over 4\alpha}\right)^{\log_7(2)\cdot{\bf p}(E)}~\leq~\mathcal{N}_{\alpha}\left(B_\rho(z,R)~\Big|~E\right)~\leq~\left({2R\over\alpha}\right)^{{\bf d}(E)}
\eeq
and 
\bel{ese4}
\left({R\over 2\alpha}\right)^{ \log_7(2)\cdot{\bf p}(E)}~\leq~\mathcal{M}_{\alpha}\left(B_\rho(z,R)~\Big|~E\right)~\leq~\left({4R\over\alpha}\right)^{{\bf d}(E)}.
\eeq
\subsection{Functions of bounded total generalized  variation} In this subsection, we  now introduce the concept of total generalized  variation of the function $g:[a,b]\to E$ which was well-studied in  \cite{MO} for the case $E=\mathbb{R}$. Consider a convex function $\Psi: [0,+\infty)\to [0,+\infty)$ such that 
\bel{Psi-cond}
\Psi(0)~=~0\qquad\mathrm{and}\qquad\Psi(s)~>~0\qquad \forall s>0\,.
\eeq
\begin{definition} The $\Psi$-total variation of $g$ over $[a,b]$ is defined as
\bel{phi-var}
TV^{\Psi}\left(g, [a,b]\right)~=~\sup_{n\in\mathbb{N}, a=x_0<x_1<...<x_n=b}\sum_{i=0}^{n-1}\Psi \left(\rho(g(x_{i}),g(x_{i+1}))\right).
\eeq
If the supremum is finite then we say that $g$ has bounded $\Psi$-total variation and denote it by $g\in BV^{\Psi}([a,b],E)$. In the case of  $\Psi(x)=|x|^{\gamma}$ for some $\gamma\geq 1$, we shall denote by 
$$
BV^{{1\over \gamma}}([a,b],E):= BV^{\Psi}([a,b]),\qquad TV^{{1\over \gamma}}\left(g, [a,b]\right)~:=~TV^{\Psi}\left(g, [a,b]\right)
$$ 
the fractional BV space on $[a,b]$ and the $\gamma$-total variation of $g$, respectively.
\end{definition}
For any function $g\in BV^{\Psi}([a,b],E)$, it is easy to show by a contradiction argument that $g$ is a regulated function, i.e., the left and right hand side limits of $g$ at $x_0\in [a,b]$ always exist, denoted by
\begin{equation*}
g(x_0-)~:=~\lim_{x\to x_0-}~g(x)\qquad\mathrm{and}\qquad g(x_0+)~:=~\lim_{x\to x_0+}~g(x).
\end{equation*}
Moreover, the set of discontinuities of $g$ 
\[
\mathcal{D}_{g}~:=~\left\{x\in [a,b]~\big|~g(x+)=g(x)=g(x-)~~\text{does not hold}\right\}
\]
is at most countable. In particular, one has the following:
\medskip
\begin{lemma}\label{repre-g} For any function $g\in BV^{\Psi}([a,b],E)$, the following function
\[
\tilde{g}(b)~=~g(b),\qquad \tilde{g}(x)~:=~g(x+)\qquad\forall x\in [a,b)
\] 
is a continuous function from the right on the interval $[a,b)$ and belongs to $BV^{\Psi}([a,b],E)$ with
\bel{TV-L1}
\rho_{{\bf L}^1}(\tilde{g},g)~=~0\qquad\mathrm{and}\qquad TV^{\Psi}\left(\tilde{g}, [a,b]\right)~\leq~TV^{\Psi}\left(g, [a,b]\right).
\eeq
\end{lemma}
{\bf Proof.} Since $\mathcal{D}_g$ is at most countable, it holds that  
\[
\rho_{{\bf L}^1}(\tilde{g},g)~=~\int_{[a,b]\backslash \mathcal{D}_g}\rho(\tilde{g}(x),g(x))dx~=~0~.
\]
On the other hand, for any partition $\{a=x_0<x_1<\dots<x_n=b\}$ of $[a,b]$,
\begin{equation*}
\sum_{i=0}^{n-1}\Psi(\rho(\tilde{g}(x_{i+1}),\tilde{g}(x_i)))~=~\Psi(\rho(g(b),g(x_{n-1}+)))+\sum_{i=0}^{n-2}\Psi(\rho(g(x_{i+1}+),g(x_i+)))~\leq~TV^{\Psi}(g,[a,b])
\end{equation*}
and this yields the second inequality in (\ref{TV-L1}).
\endproof 
\medskip

The following remark is used in the proof of  the upper estimate in Theorem \ref{main}.
\begin{remark}\label{rm1}
 Under the assumption (\ref{Psi-cond}), the function $\Psi$ is strictly increasing on $[0,+\infty)$ and 
 \bel{rtt}
 \Psi(s)~\leq~{s\over t}\cdot \Psi(t)\quad\forall 0\leq s< t~.
\eeq
 Moreover, its inverse $\Psi^{-1}$  is also strictly increasing, concave and the map $s\longmapsto\ds{\Psi^{-1}(s)\over s}$ is strictly decreasing on $[0,+\infty)$.
\end{remark}
{\bf Proof.} By the convexity of $\Psi$ and (\ref{Psi-cond}),
\[
\Psi(s)~\leq~{t-s\over t}\cdot \Psi(0)+{s\over t}\cdot \Psi(t)~=~{s\over t}\cdot \Psi(t)~<~\Psi(t)
\]
for all $0\leq s<t$. Thus, $\Psi$ is strictly increasing and convex in $[0,+\infty)$ and this implies that its inverse $\Psi^{-1}$ exists, is strictly  increasing and concave. In particular, 
\begin{align*}
\frac{\Psi^{-1}(s)}{s}~=~\frac{\Psi^{-1}(s)-\Psi^{-1}(0)}{s}~>~\frac{\Psi^{-1}(r)}{r}\qquad\forall 0< s<r
\end{align*}
and this yields the decreasing property of  the map $s\longmapsto\ds{\Psi^{-1}(s)\over s}$~.
\endproof
\section{The $\ve$-entropy for a class of $BV^{\Psi}$ functions} 
\subsection{Main results}
\setcounter{equation}{0}
Throughout this subsection, the metric space $(E,\rho)$ is assumed to be totally bounded.  For convenience, we use the notation
\[
{\bf H}_\alpha~:=~\log_2{\bf N}_\alpha\quad\mathrm{and}\qquad {\bf K}_{\alpha}~:=~\log_2{\bf M}_{\alpha}
\]
where  ${\bf N}_\alpha~:=~\mathcal{N}_\alpha(E|E)$ and ${\bf M}_{\alpha}:=\mathcal{M}_{\alpha}(E|E)$ are the $\alpha$-covering  and the $\alpha$-packing numbers of $E$ in $(E,\rho)$ and 
\[
\begin{cases}
{\bf d}&:=~{\bf d}(E)~~\mathrm{the~doubling~dimension~of~}E,\cr
{\bf p}&:=~{\bf p}(E)~~\mathrm{the~packing~dimension~of~}E.
\end{cases}
\]
Given two constants $L,V>0$, we shall establish both upper and lower  estimates on the $\ve$-entropy  of a class of uniformly bounded $\Psi$-total variation functions defined on $[0,L]$ and taking values in $(E,\rho)$,
\bel{F}
\mathcal{F}^{\Psi}_{[L,V]}~:=~\left\{f\in BV^{\Psi}\left([0,L],E\right)~\big|~ TV^{\Psi}(f,[0,L])\leq V\right\},
\eeq
in ${\bf L}^1([0,L],E)$. 
\medskip

\begin{theorem}\label{main} Assume that the function $\Psi:[0,+\infty)\to [0,+\infty)$  is convex and satisfies the condition \eqref{Psi-cond}. Then, for every $0<\ve \leq 2L\Psi^{-1}\left({V\over 4}\right)$, it holds
\begin{equation}\label{m1}
{{\bf p}V\over 2\log_2(7)\cdot \Psi\left(256\ve\over L\right)}+{\bf K}_{258\ve\over L}~\leq~ \mathcal{H}_{\ve}\left(\mathcal{F}^{\Psi}_{[L,V]}\Big|{\bf L}^1([0,L],E)\right) ~\leq~ \left[3{\bf{d}}+\log_2(5e)\right]\cdot \frac{2V}{\Psi\left({\ve\over 2L}\right)}+{\bf H}_{{\ve\over 4L}}~.
\end{equation}
 \end{theorem}
 
As a consequence, the minimal number of functions needed to represent a function in $\mathcal{F}^{\Psi}_{[L,V]}$ up to an accuracy $\ve$ with respect to ${\bf L}^1$-distance is of the order $\ds{1\over \Psi(O(\ve))}$. Indeed, from (\ref{ese3}) and (\ref{ese4}), it holds that
\[
\begin{cases}
{\bf H}_{\alpha}&\leq~{\bf d}\cdot \log_2\left(\ds\mathrm{diam}(E)\cdot {2\over \alpha}\right)\cr\cr
{\bf K}_{\alpha}&\geq~\ds{\bf p}\cdot (\log_72)\cdot \log_2\left(\ds\mathrm{diam}(E)\cdot {1\over 2\alpha}\right)
\end{cases}
\qquad\qquad\forall \alpha>0,
\]
and  (\ref{m1}) implies 
\begin{multline}\label{est1111}
{{\bf p}V\over 2\log_2(7)\cdot \Psi\left({256\ve\over L}\right)}+{\bf p}\cdot\log_7\left(\mathrm{diam}(E)\cdot {L\over 516 \ve}\right)~\leq~ \mathcal{H}_{\ve}\left(\mathcal{F}^{\Psi}_{[L,V]}~\Big|~{\bf L}^1([0,L],E)\right) \\~\leq~\left[3{\bf{d}}+\log_2(5e)\right]\frac{2V}{\Psi\left({\ve\over 2L}\right)} +{\bf d}\cdot \log_2\left({\mathrm{diam(E)}\cdot {8L\over \ve}}\right).
\end{multline}

On the other hand, one also obtains a sharp estimate on the $\ve$-entropy  for a class of uniformly bounded $\gamma$-total variation functions, i.e. $\Psi(x) = |x|^\gamma$, for all $\gamma\geq 1$. More precisely, let us denote by 
\bel{F-gamma}
\mathcal{F}^{\gamma}_{[L,V]}~=~\left\{f\in BV^{{1\over \gamma}}([0,L],E)~\big|~ TV^{{1\over \gamma}}(f,[0,L])\leq V\right\},
\eeq
it follows directly from Theorem \ref{main} that
\medskip

\begin{corollary}\label{ccc1}
For every  $0<\ve\leq\ds 2^{{\gamma -2}\over \gamma}LV^{1\over \gamma}$,  
\begin{multline}\label{m2}
{{\bf p}\over 2^{8\gamma+1}\log_2(7)}\cdot {L^{\gamma}V\over \ve^{\gamma}}+{\bf p}\cdot\log_7\left(\mathrm{diam}(E)\cdot {L\over 516\ve}\right)~~\leq~~ \mathcal{H}_{\ve}\left(\mathcal{F}^{\gamma}_{[L,V]}~|~{\bf L}^1([0,L],E)\right)\\
~\leq~2^{\gamma+1}\cdot\left[3{\bf{d}}+\log_2(5e)\right]\frac{L^\gamma V}{\ve^\gamma}+{\bf d}\cdot \log_2\left({\mathrm{diam(E)}\cdot {8L\over \ve}}\right).
\end{multline}
\end{corollary}
In particular, as $\ve$ tends to $0+$, one derives that
\begin{multline*}
{{\bf p}\over 2^{8\gamma+1}\log_2(7)}~\leq~\liminf_{\ve\to 0+}\left[{\ve^{\gamma}\over L^{\gamma}V }\cdot \mathcal{H}_{\ve}\left(\mathcal{F}^{\gamma}_{[L,V]}~|~{\bf L}^1([0,L],E)\right)\right]\\
~\leq~\limsup_{\ve\to 0+}\left[{\ve^{\gamma}\over L^{\gamma}V }\cdot \mathcal{H}_{\ve}\left(\mathcal{F}^{\gamma}_{[L,V]}~|~{\bf L}^1([0,L],E)\right)\right]~\leq~ 2^{\gamma+1} \left[3{\bf d}+\log_2(5e)\right].
\end{multline*}
Thus, the $\ve$-entropy of $\mathcal{F}^{\gamma}_{[L,V]}$  in  ${\bf L}^1([0,L],E)$ is of the order  $\ve^{-\gamma}$.
\medskip

Finally, in order to  apply our result to study the $\ve$-entropy for entropy admissible weak  solution sets  to scalar conservation laws in one-dimensional space with  weakly genuinely nonlinear fluxes, we consider the case  where the metric space $(E,\rho)$ is generated by a finite dimensional normed space $(\R^d,\|\cdot\|)$, i.e., 
\[
E~=~\R^d\qquad\mathrm{and}\qquad \rho(x,y)~=~\|x-y\|\qquad\forall x,y\in\R^d.
\]
Given an additional constant  $M>0$, the following provides upper and lower estimates for the $\ve$-entropy of  a class of uniformly bounded $\Psi$-total variation functions taking values in the open ball $B^{d}(0,M)\subset \R^d$,
\bel{F3}
\mathcal{F}^{\Psi}_{[L,M,V]}~:=~\left\{f\in BV^{\Psi}\left([0,L],B^d(0,M)\right)~\big|~ TV^{\Psi}(f,[0,L])\leq V\right\},
\eeq
in the normed space ${\bf L}^1(\R^d)$.
\medskip

\begin{corollary}\label{Eucli-cases}
Under the same assumptions in Theorem \ref{main}, it holds 
\begin{multline}\label{m3}
{Vd\over 2\log_2(7)\cdot\Psi\left({256\ve\over L}\right)}+d\cdot \log_7\left({LM\over 258\ve}\right)~\leq~ \mathcal{H}_{\ve}\left(\mathcal{F}^{\Psi}_{[L,M,V]}~\Big|~{\bf L}^1([0,L],\mathbb{R}^d)\right) \\~\leq~ \left[3d\log_25+\log_2(5e)\right]\cdot {2V\over \Psi\left({\ve\over 2L}\right)}+d\cdot \log_2\left({8LM\over \ve}+1\right)
\end{multline}
for every $0<\ve \leq 2L\Psi^{-1}\left({V\over 4}\right)$.
\end{corollary}
{\bf Proof.} It is well-known (see e.g in \cite{KT}) that 
\[
d\cdot \log_2\left({r\over \alpha}\right)~\leq~\mathcal{H}_{\alpha}\left(B^d(0,r)\Big|\mathbb{R}^d\right)~\leq~ d\cdot \log_2\left({2r\over \alpha}+1\right)
\]
for any $\alpha>0$ and open ball $B^d(0,r)\subset \R^d$. In particular, 
recalling that 
 $${\bf H}_{\alpha} = \log_2\mathcal{N}_{\alpha}\left(B^d(0,M)\Big|\mathbb{R}^d\right)\qquad\mathrm{and}\qquad {\bf K}_{\alpha} = \log_2\mathcal{M}_{\alpha}\left(B^d(0,M)\Big|\mathbb{R}^d\right),$$
 we have
\[
{\bf H}_{\alpha}~\leq~d\cdot \log_2\left({2M\over \alpha}+1\right), \qquad{\bf K}_{\alpha}~\geq~{\bf H}_{\alpha}~\geq~d\cdot \log_2\left({M\over \alpha}\right),
\]
and from Definition \ref{pc-dimension}, it holds that
\[
d~\leq~{\bf p}\big(\R^d\big)~\leq~ {\bf d}\big(\mathbb{R}^d\big)~\leq~d\cdot\log_25.
\]

Using the above estimates in (\ref{m1}), one obtains (\ref{m3}).
\qed

\quad In the next two subsections, we will present the proof of Theorem \ref{main}.
\subsubsection{Upper estimate} Towards the proof of the upper bound on $\mathcal{H}_{\ve}\left(\mathcal{F}^{\Psi}_{[L,V]}~\Big|~{\bf L}^1([0,L],E)\right)$ in Theorem \ref{main}, let us extend a result on the $\ve$-entropy for a class of bounded total variation real-valued functions in the scalar case \cite{BKP} or in \cite[Lemma 2.3]{DN}. In order to obtain a  sharp upper bound, one needs to utilize the doubling dimension of the metric space $E$ and go beyond the particular cases in \cite{BKP, DN} to estimate the $\ve$-entropy for a more general case in $E$. More precisely, considering a set of bounded total variation functions taking values in $E$, which we denote by
\bel{B}
\mathcal{F}_{[L,V]}~=~\left\{f\in BV\left([0,L],E\right)~\Big|~ TV(f,[0,L])\leq V\right\},
\eeq 
the following holds.
\v
\begin{proposition}\label{sp-case}
For every $\ds 0<\ve\leq \frac{LV}{2}$ sufficiently small, it holds that
\[
\mathcal{H}_{\ve}\left(\mathcal{F}_{[L,V]}~\Big|~{\bf L}^1([0,L],E)\right)~\leq~\left[3{\bf d}+\log_2(5e)\right]\cdot{2LV\over \ve}+{\bf H}_{{\ve\over 2L}}.
\]
\end{proposition}
{\bf Proof.} The proof is divided into four steps:
\medskip

{\bf 1.} Given  two constants $N_1\in\mathbb{Z}^+$ and $h_2>0$, let us
\begin{itemize}
\item divide $[0,L]$ into $N_1$ small intervals $I_i$ with length $h_1:=\ds{L\over N_1}$ such that $I_{N_1-1}=\ds \left[(N_1-1)h_1, L\right]$ and 
\[
I_{i}~=~[ih_1, (i+1)h_1\big)\quad\forall i\in \overline{0,N_1-2}~;
\]
\item pick an optimal $h_2$-covering $A=\left\{a_1,a_2,\dots, a_{{\bf N}_{h_2}}\right\}$ of $E$, i.e.
\[
E~\subseteq~\bigcup_{i=1}^{{\bf N}_{h_2}}B_\rho(a_i,h_2)~,
\]
where ${\bf N}_{h_2}$ is the $h_2$-covering number of $E$ (see Definition \ref{pc}). 
\end{itemize}
A function $f\in \mathcal{F}_{[L,V]}$ can be approximated by a piecewise constant function $f^{\sharp}: [0,L]\to A$ defined as follows:
\[
f^{\sharp}(s)~=~a_{f,i}\quad\forall s\in I_i,~ i\in  \overline{0,N_1-1 }
\]
for some $a_{f,i}\in A$ such that $\ds f\left(t_i\right)\in B_\rho(a_{f,i},h_2)$ with $t_i:=\ds{2i+1\over 2}h_1$. Notice that $a_{f,i}$ is not a unique choice. With this construction, the ${\bf L}^1$-distance between $f$ and $f^{\sharp}$ can be bounded above by 
\begin{align*}
\rho_{\textbf{L}^1}(f, f^{\sharp})~&\leq~\sum_{i=0}^{N_1-1}\int_{I_i} \rho(f(s),f^{\sharp}\left(s\right))ds~=~\sum_{i=0}^{N_1-1}\int_{I_i} \rho(f(s),a_{f,i})ds \\
&~\leq~\sum_{i=0}^{N_1-1}\int_{I_i} \Big[\rho(f(s),f(t_i))+\rho(f(t_i),a_{f,i})\Big]ds~<~\sum_{i=0}^{N_1-1}\int_{I_i} \Big[\rho(f(s),f(t_i))+h_2\Big]ds\\
&~\leq~\left(\sum_{i=0}^{N_1-1}{|I_i|\over 2} \cdot \left[TV(f, [ih_1,t_i])+TV(f,[t_i,(i+1)h_1])\right]\right)+Lh_2\\
&~=~{h_1\over 2}\cdot TV(f,[0,L])+Lh_2~\leq~{LV\over 2N_1}+Lh_2
\end{align*}
and the total variation of $f^{\sharp}$ over $[0,L]$ can be estimated by 
\begin{eqnarray*}
TV\left(f^{\sharp},[0,L]\right)&=&\sum_{i=0}^{N_1-2} \rho(a_{f,i},a_{f,i+1})\\
&\leq&\sum_{i=0}^{N_1-2}\Big[\rho(a_{f,i+1},f(t_{i+1}))+\rho(f(t_i),a_{f,i})+\rho\left(f(t_{i+1}\right), f\left(t_i\right))\Big]\\
&\leq& \sum_{i=0}^{N_1-2}\Big[2h_2+\rho (f\left(t_{i+1}\right),f\left(t_i\right))\Big]~\leq~2(N_1-1)\cdot h_2+V.
\end{eqnarray*}
Consider the following set of piecewise constant functions 
\begin{eqnarray*}
\mathcal{F}^{\sharp}_{[N_1,h_2]}&=&\Big\{\varphi: [0,L]\to A~~\Big|~~ \varphi(s)=\varphi(t_i)~~\forall s\in I_i, i\in\overline{0,N_1-1}\cr
\cr
& &\qquad\qquad\qquad\qquad\qquad\qquad\mathrm{and}\quad  TV(\varphi, [0,L])\leq  2(N_1-1)\cdot h_2+V\Big\}.
\end{eqnarray*}
The set $\mathcal{F}_{[L,V]}$ is covered by a finite collection of closed balls centered at $\varphi \in \mathcal{F}^{\sharp}_{[N_1,h_2]}$ of radius ${LV\over 2N_1}+Lh_2 $ in ${\bf L}^1([0,L],E)$, i.e.,
\[
\mathcal{F}_{[L,V]}~\subseteq~\bigcup_{\varphi\in \mathcal{F}^{\sharp}_{[N_1,h_2]}}~\overline{B}_{{\bf L}^1([0,L],E)}\left(\varphi,{LV\over 2N_1}+Lh_2\right)
\]
and the  Definition \ref{entropy} yields 
\bel{ee1}
\mathcal{H}_{\left[{LV\over 2N_1}+Lh_2\right]}\left(\mathcal{F}_{[L,V]}~\Big|~{\bf L}^1([0,L],E)\right)~\leq~\log_2 \mathrm{Card} \left(\mathcal{F}^{\sharp}_{[N_1,h_2]}\right).
\eeq

{\bf 2.} In order to provide an upper bound on $\mathrm{Card} \left(\mathcal{F}^{\sharp}_{[N_1,h_2]}\right)$, we introduce a discrete metric $\rho^{\sharp}:A\times A\to \mathbb{N}$ associated to $\rho$ as follows: 
	\bel{rho-s}
	\rho^{\sharp}(x,y)~:=~\begin{cases}
		0\qquad\mathrm{if}\qquad x=y,\cr\cr
		q+1\qquad\mathrm{if}\qquad \ds {\rho(x,y)\over h_2}~\in~\big(q,q+1] \qquad {\text{  for some }}q\in\mathbb{N}~,
	\end{cases}
	\eeq
	for every $x$, $y\in A$.}
Since $A$ is an optimal  $h_2$-covering of $E$, one has
\[
\mathrm{Card}\left(A\bigcap B_\rho(a,r)\right)~\leq~\mathcal{N}_{h_2}\left(B_\rho(a,r+h_2)\big|E\right)\quad\forall a\in A, r>0
\]
and  the second inequality in (\ref{ese3}) yields 
\[
\mathrm{Card}\left(A\bigcap B_\rho(a,r)\right)~\leq~\left(2\cdot \left({r\over h_2}+1\right)\right)^{\bf d}.
\]
Hence, for every $\ell \geq 1$ and $x\in A$, it holds
\begin{multline}\label{nu1}
\mathrm{Card}\left(\overline{B}_{\rho^{\sharp}}(x,\ell-1)\right)~=~\mathrm{Card}\left(\{y\in A~|~\rho^{\sharp}(x,y)\leq \ell-1\}\right)\\
~ =~\mathrm{Card}\left(A\bigcap B_\rho\left(x,(\ell - 1) h_2 \right)\right)~\leq~(2 \ell)^{\bf d}.
\end{multline}
For any given  $f^{\sharp}\in  \mathcal{F}^{\sharp}_{[N_1,h_2]}$, the following increasing step function $\varphi_{f^{\sharp}}:[0,L]\to\mathbb{N}$ defined by 
\bel{varphi}
\varphi_{f^{\sharp}}(s)~=~
\begin{cases}
0\quad\forall s \in I_0\cr\cr
\ds \sum_{\ell=0}^{i-1} \rho^{\sharp}\left(f^{\sharp}(t_{\ell}),f^{\sharp}(t_{\ell+1})\right)+i-1\quad\forall s\in I_i, i\in \overline{1, N_1-1}~
\end{cases}
\eeq
measures  the total  of jumps of $f^{\sharp}$ up to time $t_i$. 
From (\ref{rho-s}), one has
\begin{multline}\label{eeet1}
\sup_{t\in [0,L]} \left|\varphi_{f^{\sharp}}(t)\right|~\leq~\sum_{\ell=0}^{N_1-2} \rho^{\sharp}\left(f^{\sharp}(t_{\ell}),f^{\sharp}(t_{\ell+1})\right)+N_1-2\\~\leq~\sum_{\ell=0}^{N_1-2} \left({\rho(f^{\sharp}(t_{\ell}),f^{\sharp}(t_{\ell+1}))\over h_2}+1\right)+N_1-2~\leq~{TV(f^{\sharp}, [0,L])\over h_2}+2N_1-3\\
\leq  {1\over h_2}\cdot \left(2(N_1-1)\cdot h_2+V\right)+2N_1-3~=~4N_1-5+{V\over h_2}~.
\end{multline}
In particular, upon setting $\Gamma_{[N_1,h_2]}:=4N_1-4+\ds\left\lfloor {V\over h_2}\right\rfloor$, a constant depending on $N_1$ and $h_2$, the function $\varphi_{f^{\sharp}}$ in (\ref{varphi}) satisfies 
\[
\varphi_{f^{\sharp}}(s)~=~\varphi_{f^{\sharp}}(t_i)~\in~\left\{0,1,2,\dots,\Gamma_{[N_1,h_2]}-1\right\}\quad\forall s\in I_i, ~i\in \overline{0,N_1-1}~.
\] 
Thus, if we  consider the map $T: \mathcal{F}^{\sharp}_{[N_1,h_2]}\to \mathcal{B}([0,L],[0,+\infty))$ such that 
\[
T(f^{\sharp})~=~\varphi_{f^{\sharp}} \qquad\forall f^{\sharp}\in \mathcal{F}^{\sharp}_{[N_1,h_2]}~,
\]
then
\[
T\left(\mathcal{F}^{\sharp}_{[N_1,h_2]}\right)~=~\left\{\varphi_{f^{\sharp}}~\Big|~f^{\sharp}\in \mathcal{F}^{\sharp}_{[N_1,h_2]}\right\}~\subseteq~\mathcal{I}_{[N_1,h_2]}~.
\]
Here, $\mathcal{I}_{[N_1,h_2]}$ is the set of increasing step functions $\phi:[0,L]\to \left\{0,1,2,\dots,\Gamma_{[N_1,h_2]}-1\right\}$ such that 
\[
\phi(0)~=~0\quad\mathrm{and}\quad \phi(s)~=~\phi(t_i)\quad\forall i\in \overline{0,N_1-1},~ s\in I_i~.
\]

Since the cardinality of $\mathcal{I}_{[N_1,h_2]}$ is equal to $\ds {\Gamma_{[N_1,h_2]}\choose N_1-1}$, one has 
\bel{eest1}
\mathrm{Card}\left(T\left(\mathcal{F}^{\sharp}_{[N_1,h_2]}\right)\right)~\leq~\mathrm{Card}(\mathcal{I}_{[N_1,h_2]})~=~{\Gamma_{[N_1,h_2]} \choose N_1 - 1}~.
\eeq
{\bf  3.} To complete the proof, we need to establish an upper estimate on the cardinality of $T^{-1}(\varphi_{f^{\sharp}})$, the set of functions in $\mathcal{F}^{\sharp}_{[N_1,h_2]}$ that have the same total length of jumps as that of $f^{\sharp}$ at any time $t_i$. In order to do so,  for any given $f^{\sharp}\in \mathcal{F}^{\sharp}_{[N_1,h_2]}$,  we set
\[
k^{\sharp}_i~:=~\ds \rho^{\sharp}\left(f^{\sharp}(t_{i}),f^{\sharp}(t_{i+1})\right)\quad\forall i\in \overline{0,N_1-2}~.
\]
As in (\ref{eeet1}), we have 
\[
\sum_{i=0}^{N_1-2}k^{\sharp}_i~=~\ds \sum_{i=0}^{N_1-2} \rho^{\sharp}\left(f^{\sharp}(t_{i}),f^{\sharp}(t_{i+1})\right)~\leq~3(N_1-1)+{V\over h_2}
\]
and
\begin{eqnarray*}
T^{-1}(\varphi_{f^{\sharp}})&=&\left\{g\in \mathcal{F}^{\sharp}_{[N_1,h_2]}~\Big|~\rho^{\sharp}\left(g(t_{i+1}),g(t_i)\right)=k^{\sharp}_i~~\forall i\in\overline{0,N_1-2}\right\}\cr\cr
&\subseteq&\left\{g\in \mathcal{F}^{\sharp}_{[N_1,h_2]}~\Big|~g(t_{i+1})\in \overline{B}_{\rho^{\sharp}}\left(g(t_i), k^{\sharp}_i\right)\forall i\in\overline{0,N_1-2}\right\}.
\end{eqnarray*}
Observe from (\ref{nu1}) that if $g(t_{i})$ is already chosen then there are at most $(2 k^{\sharp}_i)^{\bf d}$ choices for $g(t_{i+1})$. Since we have ${\bf N}_{h_2}$ choices of the starting point $g(0)$, the cardinality of  $T^{-1}(\varphi_{f^{\sharp}})$ can be estimated as follows
 \begin{multline}\label{eest2}
\mathrm{Card}\left(T^{-1}(\varphi_{f^{\sharp}})\right)~\leq~{\bf N}_{h_2}\cdot \ds \Pi_{i=0}^{N_1-2}(2 k^{\sharp}_i)^{\bf d}~\leq~{\bf N}_{h_2}\cdot\ds\left(\ds{\sum_{i=0}^{N_1-2} 2 k^{\sharp}_i\over N_1-1}\right)^{{\bf d}(N_1-1)}\\
~\leq~{\bf N}_{h_2}\cdot \left({2 \left(3(N_1-1)+{V\over h_2}\right)\over N_1-1}\right)^{{\bf d}(N_1-1)}~=~{\bf N}_{h_2}\cdot\left(6+{2\over N_1-1}\cdot {V\over h_2}\right)^{{\bf d}(N_1-1)}.
\end{multline}
Recalling (\ref{eest1})-(\ref{eest2}) and the classical Stirling's approximation $$(N_1-1)!~\geq~\sqrt{2\pi (N_1-1)}\cdot  \left({N_1-1\over e}\right)^{N_1-1}~,$$ we estimate
\begin{multline*}
\mathrm{Card} \left(\mathcal{F}^{\sharp}_{[N_1,h_2]}\right)~\leq~{\bf N}_{h_2}\cdot\left(6+{2\over N_1-1}\cdot {V\over h_2}\right)^{{\bf d}(N_1-1)}\cdot {\Gamma_{[N_1,h_2]}\choose N_1-1}\\
~=~{\bf N}_{h_2}\cdot\left(6+{2\over N_1-1}\cdot {V\over h_2}\right)^{{\bf d}(N_1-1)}\cdot {\left(\Gamma_{[N_1,h_2]}-N_1+2\right)\dots \Gamma_{[N_1,h_2]}\over (N_1-1)!}\\
~\leq~{{\bf N}_{h_2}\over \sqrt{2\pi (N_1-1)}}\cdot\left(6+ {2\over N_1-1}\cdot {V\over h_2}\right)^{{\bf d}(N_1-1)}\cdot \left({\Gamma_{[N_1,h_2]}\over N_1-1}\right)^{N_1-1}\cdot e^{N_1-1} \\
~\leq~{\bf N}_{h_2}\cdot\left(6+ {2\over N_1-1}\cdot  {V\over h_2}\right)^{{\bf d}(N_1-1)}\cdot \left(4e+{V\over h_2}\cdot {e\over N_1-1}\right)^{(N_1-1)}.
\end{multline*}
Thus, (\ref{ee1}) yields 
\begin{multline}\label{ee2}
\mathcal{H}_{\left[{LV\over 2N_1}+{Lh_2 }\right]}\left(\mathcal{F}_{[L,V]}~\Big|~{\bf L}^1([0,L],E)\right)~\leq~{\bf d}\cdot (N_1-1)\cdot\log_2\left(6+ {V\over h_2}\cdot {2\over N_1-1}\right)\\
+(N_1-1)\cdot \log_2\left(4e+ {V\over h_2}\cdot {e\over N_1-1}\right)+ {\bf H}_{h_2}~.
\end{multline}

{\bf  4.} For every $0< \ve \leq \ds {LV\over 2}$, by choosing $N_1\in \mathbb{Z}^+$ and $h_2>0$ such that 
\[
{3LV\over 2\ve}~<~N_1-1~=~\left\lfloor{3LV\over 2\ve}\right\rfloor+1~\leq~{2LV\over \ve}~,\qquad h_2~=~{V\over N_1-1}~,
\]
we have
\[
{LV\over 2N_1}+{Lh_2 }~\leq~{LV\over 2N_1}+{LV\over N_1-1}~\leq~{3LV\over 2(N_1-1)}~<~\ve\qquad\mathrm{and}\qquad h_2~\geq~{\ve\over 2L}~.
\]
Thus, (\ref{ee2}) implies that 
\[
\mathcal{H}_{\ve}\left(\mathcal{F}_{[L,V]}~\Big|~{\bf L}^1([0,L],E)\right)~\leq~\left[3{\bf d}+\log_2(5e)\right]\cdot{2LV\over \ve}+{\bf H}_{{\ve\over 2L}}
\]
and this completes the proof.  
\endproof
\v

Using Proposition \ref{sp-case}, we now proceed to provide a proof for the upper estimate of the $\ve$-entropy for the set $\mathcal{F}^{\Psi}_{[L,V]}$ in ${\bf L}^1([0,L], E)$.
\v

\indent {\bf Proof of the upper estimate in Theorem \ref{main}.} From Lemma \ref{repre-g}, one has 
\bel{rl1}
 \mathcal{H}_{\ve}\left(\mathcal{F}^{\Psi}_{[L,V]}~\Big|~{\bf L}^1([0,L],E)\right)~=~ \mathcal{H}_{\ve}\left(\tilde{\mathcal{F}}^{\Psi}_{[L,V]}~\Big|~{\bf L}^1([0,L],E)\right)
\eeq
with $\tilde{\mathcal{F}}^\Psi_{[L,V]}=\left\{f\in\mathcal{F}^{\Psi}_{[L,V]}~\Big|~f~\mathrm{is~continuous~from~the~right~on~the~interval~} [0,L)\right\}$. Thus, it is sufficient to prove the second inequality in (\ref{m1}) for $\tilde{\mathcal{F}}^{\Psi}_{[L,V]}$ instead of $\mathcal{F}^{\Psi}_{[L,V]}$.
\medskip

{\bf 1.} For a fixed constant $h>0$ and $f\in \tilde{\mathcal{F}}^{\Psi}_{[L,V]}$, let $A_{f,h}=\left\{x_0,x_1,x_2,...,x_{N_{f,h}}\right\}$ be a partition of $[0,L]$ which is defined by induction as follows:
\begin{equation}\label{xn}
x_0~=~0,\qquad x_{i+1}~=~ \sup\left\{x\in (x_i,L)~\big|~\rho(f(y),f(x_i))\in[0,h]\quad\forall y\in  (x_i,x]\right\}
\end{equation}
for all $ i\in\overline{0,N_{f,h}-1}$. Since $f$ is continuous from the right on $[0,L)$, it holds 
\[
\rho(f(x_{i}),f(x_{i+1}))~\geq~h\qquad\forall i\in\overline{0,N_{f,h}-2}.
\]
Thus, the increasing property of $\Psi$ implies that 
\begin{equation*}
V~\geq~TV^{\Psi}(f,[0,L])~\geq~\sum_{i=0}^{N_{f,h}-2}\Psi\left(\rho(f(x_i),f(x_{i+1}))\right)~\geq~(N_{f,h}-1)\cdot\Psi(h),
\end{equation*}
and this yields 
\bel{bb}
N_{f,h}-1~\leq~{TV^{\Psi}(f,[0,L])\over \Psi(h)}~\leq~{V\over \Psi(h)}~<~+\infty~.
\eeq

Introduce a piecewise constant function $f_h: [0,L]\to E$ such that 

\[
f_h(x)~=~\begin{cases}
f(x_i)&\qquad\forall x\in [x_i,x_{i+1})~,~ i\in\overline{0,N_{f,h}-2}\cr\cr
f\left(x_{N_{f,h}-1}\right)&\qquad\forall x\in \left[x_{N_{f,h}-1},L]\right..
\end{cases}
\]

From (\ref{xn}), the ${\bf L}^1$-distance between $f_h$ and $f$ is bounded by 
\begin{multline}\label{l1-d}
\rho_{\textbf{L}^1}(f_h,f)~=~\int_{[0,L]}\rho(f_h(x),f(x))dx~=~\sum\limits_{i=0}^{N_{f,h}-1}\int_{[x_{i},x_{i+1})}\rho(f(x_i),f(x))dx\\
~\leq~h\cdot \sum\limits_{i=0}^{N_{f,h}-1}(x_{i+1}-x_{i})~=~Lh\,.
\end{multline}

On the other hand, by the convexity of $\Psi$ we have 
\begin{eqnarray*}
V~\geq~\sum_{i=0}^{N_{f,h}-2}\Psi\left(\rho(f(x_i),f(x_{i+1}))\right)&\geq&\left(N_{f,h}-1\right)\cdot \Psi\left({1\over N_{f,h}-1}\cdot \sum_{i=0}^{N_{f,h}-2}\rho(f(x_i),f(x_{i+1})) \right)\cr\cr
&=&\left(N_{f,h}-1\right)\cdot\Psi\left({TV(f_h,[0,L))\over N_{f,h}-1}\right)
\end{eqnarray*}
and the strictly increasing property of $\Psi^{-1}$ implies 
\[
TV(f_h,[0,L))~\leq~(N_{f,h}-1)\cdot \Psi^{-1}\left({V\over N_{f,h}-1}\right).
\]
From Remark \ref{rm1} and (\ref{bb}), it holds that
\[
\Psi^{-1}\left(V\over N_{f,h}-1\right)\cdot {N_{f,h}-1\over V}~\leq~\Psi^{-1}\left(\Psi(h)\right)\cdot{1\over \Psi(h)}~=~{h\over\Psi(h)}
\]
and this yields 
\[
TV(f_h, [0,L])~\leq~{h\over \Psi(h)}\cdot V~=:~V_h~.
\]
From (\ref{l1-d}) and (\ref{B}), the set  $\tilde{\mathcal{F}}^{\Psi}_{[L,V]}$ is covered by a collection of closed balls centered at $g \in \mathcal{F}_{[L,V_h]}$ of radius $Lh$ in ${\bf L}^1([0,L],E)$, i.e.,
\[
\tilde{\mathcal{F}}^{\Psi}_{[L,V]}~\subseteq~\bigcup_{g\in\mathcal{F}_{[L,V_h]} }\overline{B}_{{\bf L}^1([0,L],E)}(g,Lh).
\]

In particular, for every $\ve>0$, choosing $h={\ve\over 2L}$ we have 
\[
V_{\ve\over 2L}~=~{\ve V\over 2L\cdot \Psi\left({\ve\over 2L}\right)}\qquad\mathrm{and}\qquad \tilde{\mathcal{F}}^{\Psi}_{[L,V]}~\subseteq~\bigcup_{g\in\mathcal{F}_{\left[L,V_{\ve\over 2L}\right]} }\overline{B}_{{\bf L}^1([0,L],E)}\left(g,{\ve\over 2}\right)
\]
and this implies
\bel{ee3}
 \mathcal{H}_{\ve}\left(\mathcal{\tilde{F}}^{\Psi}_{[L,V]}~\Big|~{\bf L}^1([0,L],E)\right)~\leq~\mathcal{H}_{{\ve\over 2}}\left(\mathcal{F}_{\left[L,V_{\ve\over 2L}\right]}~\Big|~{\bf L}^1([0,L],E)\right).
\eeq
If $0<\ds\ve \leq 2L\Psi^{-1}\left({V\over 4}\right)$ then 
\[
\ve~\leq~\ve\cdot {V\over 4\cdot \ds\Psi\left({\ve\over 2L}\right)}~=~{L\over2}\cdot  {\ve V\over 2L\cdot \ds\Psi\left({\ve\over 2L}\right)}~=~{L\over2}\cdot V_{\ve\over 2L}~.
\]

In this case, one can apply Proposition \ref{sp-case} to get
\begin{eqnarray*}
\mathcal{H}_{\frac{\ve}{2}}\left({\mathcal{F}}_{\left[L,V_{\ve\over 2L}\right]}~\Big|~{\bf L}^1([0,L],E)\right)&\leq&\left[3{\bf{d}}+\log_2(5e)\right]\cdot \frac{4LV_{\ve\over 2L}}{\ve}+{\bf H}_{{\ve\over 4L}}\\
&=&\left[3{\bf{d}}+\log_2(5e)\right]\cdot \frac{2V}{\Psi\left({\ve\over 2L}\right)}+{\bf H}_{{\ve\over 4L}}
\end{eqnarray*}
and thereafter, we use (\ref{rl1}), (\ref{ee3}) to obtain the second inequality in (\ref{m1}).
\endproof
\subsubsection{Lower estimate} To prove the first inequality  in Theorem \ref{main}, let us provide a lower estimate on the $\ve$-entropy in ${\bf L}^1([0,L],E)$ to   
\bel{ddG}
\mathcal{G}^{\Psi}_{[L,V,h,x]}:=\left\{g:[0,L]\to B_\rho(x,h)~\Big|~TV^{\Psi}(g, [0,L])\leq V\right\},
\eeq
a class of bounded $\Psi$-total variation functions over $[0,L]$ taking values in the ball centered at a point $x \in E$ of radius $h>0$. 
\v

\begin{lemma} Assume that ${\bf p}\geq 1$.  For every $\ve>0$, it holds 
\bel{low-11}
\mathcal{M}_{\ve}\left(\mathcal{G}^{\Psi}_{\left[L,V,2^{(4+2/{\bf \tilde{p}})}\cdot{\ve\over L},x\right]}~\Big|~{\bf L}^1([0,L],E)\right)~\geq~ 2^{{\bf \tilde{p}}V\over 2 \Psi\left(2^{(4+2/{\bf \tilde{p}})}\cdot {2\ve\over L}\right)}
\eeq
where ${\bf \tilde{p}}=\log_7(2)\cdot {\bf p}$.
\end{lemma}
{\bf Proof.} The proof is divided into two steps:
\medskip

{\bf  1.} We first recall from (\ref{ese4}) that 
\[
{\mathcal{M}_{2^{-(2+2/{\bf \tilde{p}})}\cdot h}(B_{\rho}(x,h)|E)}~\geq~ \left({h\over 2\cdot 2^{-(2+2/{\bf \tilde{p}})}\cdot h}\right)^{\bf {\tilde{p}}}~=~2^{{\bf \tilde{p}}+2}\qquad\forall h>0~.
\] 
Given two constants $h>0$ and $N_1\in \mathbb{Z}^+$, let us
\begin{itemize} 
\item divide $[0,L]$ into $N_1$ small mutually disjoint intervals $I_i$ with length $\ds h_1={L\over N_1}$ as in Proposition \ref{sp-case};
\item take a  $\ds \left(2^{-(2+2/{\bf \tilde{p}})}\cdot h\right)-$packing $A_h=\left\{a_1,a_2,\dots, a_{2^{\bf \tilde{p}+2}}\right\}$ of $B_{\rho}(x,h)$, i.e., 
\[
A_h~\subseteq~B_{\rho}(x,h)\qquad\mathrm{and}\qquad \rho(a_i,a_j)~>~2^{-(2+2/{\bf \tilde{p}})}\cdot h
\]
for all $ a_i\neq a_j \in A_h$~.
 \end{itemize}
Consider the set of indices  
$$
\Delta_{h,N_{1}} = \Big\lbrace \delta = (\delta_i)_{i \in \lbrace 0,1, ~\cdots,~N_{1} -1\rbrace}~\Big|~\delta_i \in A_h \Big\rbrace
$$
and define a class of piecewise constant functions on $[0,L]$ as follows:
\[
\mathcal{G}_{h,N_{1}}~=~\left\{g_{\delta}=\sum_{i=0}^{N_1-1}\delta_i \cdot \chi_{I_i}~\Big|~\delta\in \Delta_{h,N_{1}}\right\}.
\]
For any  $\delta\in \Delta_{h,N_{1}}$, the $\Psi$-total variation of $g_{\delta}$ is bounded by  
\[
TV^{\Psi}(g_{\delta}, [0,L])~\leq~(N_1-1)\cdot \Psi(2h).
\]
%
Hence, under the following condition on $h$ and $V$
\bel{Constra}
(N_1-1)\cdot \Psi(2h)~\leq~V, 
\eeq
the definition of $\ds\mathcal{G}^{\Psi}_{[L,V,h,x]}$ in (\ref{ddG}) implies that $g_{\delta}\in \ds\mathcal{G}^{\Psi}_{[L,V,h,x]}$ for every $\delta \in \Delta_{h,N_{1}}$ and thus 
$$\mathcal{G}_{h,N_{1}}~ \subseteq~\ds\mathcal{G}^{\Psi}_{[L,V,h,x]}~.$$
%
In particular, we get 
\bel{ineq1}
\mathcal{M}_{\varepsilon}\Big(\mathcal{G}^{\Psi}_{[L,V,h,x]}~\big|~ {\bf{L^1}}([0,L],E) \Big)~\geq~ \mathcal{M}_{\varepsilon}\Big(\mathcal{G}_{h,N_{1}}~\big|~ {\bf{L^1}}([0,L],E) \Big)\quad \text{ for all } \varepsilon > 0~.
\eeq

{\bf  2.} Let us  provide a lower bound on the $\ve$-packing number  $\mathcal{M}_{\varepsilon}\Big(\mathcal{G}_{h,N_{1}}~\big|~ {\bf{L^1}}([0,L],E) \Big)$. For any given $\delta,\tilde{\delta}\in  \Delta_{h,N_{1}}$ and $\ve>0$, we define 
\[
\mathcal{I}_{\tilde{\delta}} (2\varepsilon)~=~\Big\lbrace \delta \in \Delta_{h,N_{1}}~\Big|~\rho_{\bf{L}^1}(g_\delta, g_{\tilde{\delta}}) \leq 2\varepsilon \Big\rbrace,\quad \eta(\delta, \tilde{\delta})~=~\text{Card}\left(\left\lbrace i \in \overline{0,N_1-1}~\big|~ \delta_i \neq \tilde{\delta}_i \right\rbrace\right).
\]
The ${\bf L}^1$-distance between $g_{\delta}$ and $g_{\tilde{\delta}}$ is bounded below by
\begin{eqnarray*}
\rho_{\bf{L}^1}(g_{\delta},g_{\tilde{\delta}})&=&\sum_{i=0} ^{N_1-1}\int_{I_i}\rho(g_{\delta}(t),g_{\tilde{\delta}}(t))\,dt~=~~\sum_{i=0}^{N_1-1}\rho(\delta_i,\tilde{\delta}_i)\cdot |I_i|\cr\cr
&=&{L\over N_1}\cdot \sum_{i=0}^{N_1-1}\rho(\delta_i,\tilde{\delta}_i)~>~2^{-(2+2/{\bf \tilde{p}})}\cdot {Lh\over N_1}\cdot \eta(\delta,\tilde{\delta})
\end{eqnarray*}
and this implies the  inclusion
\bel{daad}
\mathcal{I}_{\tilde{\delta}} (2\varepsilon)~\subseteq~\left\lbrace \delta \in \Delta_{h,N_{1}}~\Big|~ \eta(\delta,\tilde{\delta}) < \frac{2^{3+2/{\bf \tilde{p}}}N_1\varepsilon}{Lh} \right\rbrace.
\eeq
On the other hand, for every $r\in \overline{0,N_1-1}$, we compute
\[
\mathrm{Card}\left(\Big\lbrace \delta \in \Delta_{h,N_{1}}~\Big|~ \eta(\delta,\tilde{\delta})=r\Big\rbrace\right)~=~{{N_{1}} \choose {r}}\cdot \left(2^{{\bf \tilde{p}}+2}-1\right)^{r}.
\]
Thus, (\ref{daad}) implies that 
\[
\mathrm{Card}\left(\mathcal{I}_{\tilde{\delta}} (2\varepsilon)\right)~\leq~\mathrm{Card}\left(\Big\lbrace \delta \in \Delta_{h,N_{1}}~\Big|~ \eta(\delta,\tilde{\delta})<\frac{2^{3+2/{\bf \tilde{p}}}N_1\varepsilon}{Lh} \Big\rbrace\right)~\leq~\sum\limits _{r=0}^{\left \lfloor \frac{2^{3+2/{\bf \tilde{p}}}N_1\varepsilon}{Lh}\right \rfloor}{{N_{1}} \choose {r}}\cdot \left(2^{{\bf \tilde{p}}+2}-1\right)^{r}.
\]
In particular, for every $0<\ve\leq \ds 2^{-(4+2/{\bf \tilde{p}})}Lh$, we have 
\begin{align}\label{te2}
\mathrm{Card}\left(\mathcal{I}_{\tilde{\delta}} (2\varepsilon)\right)~\leq~\sum\limits _{r=0}^{\left \lfloor \frac{N_1}{2}\right \rfloor}{{N_{1}} \choose {r}}\cdot \left(2^{{\bf \tilde{p}}+2}-1\right)^{r}~\leq~\left(2^{{\bf \tilde{p}}+2}-1\right)^{\frac{N_1}{2}}\cdot \sum\limits _{r=0}^{\left \lfloor \frac{N_1}{2}\right \rfloor}{{N_{1}} \choose {r}} \notag \\
~\leq~2^{({{\bf \tilde{p}}+2})\frac{N_1}{2}}\cdot2^{N_1}~=~ 2^{N_1 (2+{\bf \tilde{p}}/2)}~.
\end{align}

\medskip
Recalling Definition \ref{pc}, we then obtain that 
\[
\mathcal{M}_{\varepsilon}\Big(\mathcal{G}_{h,N_{1}}~\big|~ {\bf{L^1}}([0,L],E) \Big)~\geq~{\mathrm{Card}\left(\mathcal{G}_{h,N_{1}}\right)\over \mathrm{Card}\left(\mathcal{I}_{\tilde{\delta}} (2\varepsilon)\right)}~\geq~{2^{N_1({\bf \tilde{p}}+2)}\over 2^{N_1(2+{\bf \tilde{p}}/2)}}~=~2^{N_1{\bf \tilde{p}}/2}~.
\]

Finally, by choosing $h=\ds 2^{(4+2/{\bf \tilde{p}})}\cdot {\ve\over L}$ and $N_1=\ds\left\lfloor{V\over \Psi(2^{(4+2/{\bf \tilde{p}})}\cdot {2\ve\over L})}\right\rfloor+1$ such that (\ref{Constra}) holds, we derive 
\[
\mathcal{M}_{\varepsilon}\Big(\mathcal{G}_{2^{(4+2/{\bf \tilde{p}})}\cdot {\ve\over L},N_{1}}~\big|~ {\bf{L^1}}([0,L],E) \Big)~\geq~\ds 2^{{\bf \tilde{p}}V\over 2\Psi\left(2^{(4+2/{\bf \tilde{p}})}\cdot {2\ve\over L}\right)}
\] 
and thereafter, (\ref{ineq1}) yields (\ref{low-11}).
\endproof
\quad\\
To complete this section, we prove the first inequality in (\ref{m1}).
\newpage
{\bf Proof of the lower bound in Theorem \ref{main}.} For any $0<2h<h_2$,
%
let  $\{x_1,x_2,\dots, x_{{\bf M}_{h_2}}\}\subseteq E$ be an $h_2$-packing of $E$ with size ${\bf M}_{h_2}$, i.e.,  
\[
B_{\rho}\left(x_i,{h_2\over 2}\right)\bigcap B_{\rho}\left(x_j,{h_2\over 2}\right)~=~\emptyset\quad\forall i\neq j \in \overline{1,{\bf M}_{h_2}}~.
\]
Recalling the definition of $\mathcal{G}^{\Psi}_{[L,V,h,x]}$ in (\ref{ddG}), we have 
\begin{eqnarray*}
\rho_{\bf{L}^1}(f_i,f_j)&\geq&\int_{[0,L]}\Big[\rho(x_i,x_j)-\rho(x_i,f_i(s))-\rho(x_j,f_j(s))\Big]\,ds~\geq~L\cdot (h_2-2h)=:L_{h,h_2}
\end{eqnarray*}
for any $f_i\in \mathcal{G}^{\Psi}_{[L,V,h,x_i]}$ and $f_j\in \mathcal{G}^{\Psi}_{[L,V,h,x_j]}$ with $i\neq j\in\overline{1,{\bf M}_{h_2}}$~. Thus, Lemma \ref{Rl} implies that  
\begin{eqnarray*}
\mathcal{N}_{{L_{h,h_2}\over 2}}\left(\mathcal{F}^{\Psi}_{[L,V]}~\Big|~{\bf L}^1([0,L],E)\right)&\geq&\mathcal{M}_{L_{h,h_2}}\left(\mathcal{F}^{\Psi}_{[L,V]}~\Big|~{\bf L}^1([0,L],E)\right)\cr\cr
&\geq&\mathcal{M}_{L_{h,h_2}}\left(\bigcup_{i=1}^{{\bf M}_{h_2}} \mathcal{G}^{\Psi}_{[L,V,h,x_i]}~\Big|~{\bf L}^1([0,L],E)\right)\cr\cr
&=&\sum_{i=1}^{{\bf M}_{h_2}}\mathcal{M}_{L_{h,h_2}}\left(\mathcal{G}^{\Psi}_{[L,V,h,x_i]}~\Big|~{\bf L}^1([0,L],E)\right).
\end{eqnarray*}
Two cases are considered:

$\bullet$ If ${\bf p}=0$ then by choosing $h=\ds {\ve\over L}$ and $h_2=\ds {4\ve\over L}$ such that $L_{h,h_2}=2\ve$, we have 
\[
\mathcal{N}_{\ve}\left(\mathcal{F}^{\Psi}_{[L,V]}~\Big|~{\bf L}^1([0,L],E)\right)~\geq~{\bf M}_{4\ve\over L}
\]
and this particularly implies the first inequality in (\ref{m1}).

$\bullet$ Otherwise if ${\bf p}\geq 1$, then for any $\ve>0$, choosing $\ds h= 2^{(5+2/{\bf \tilde{p}})}\cdot {\ve\over L}$ and $\ds h_2=\left(2+2^{\left(6+2/{\bf \tilde{p}}\right)}\right)\cdot {\ve\over L}$ with ${\bf \tilde{p}}=\ds\log_7(2)\cdot {\bf p}$ such that $L_{h,h_2}=2\ve$, we can apply (\ref{low-11}) to $\mathcal{G}^{\Psi}_{[L,V,h,x_i]}$ for every $i\in\overline{1,{\bf M}_{h_2}}$ to obtain
\begin{multline*}
\mathcal{N}_{\ve}\left(\mathcal{F}^{\Psi}_{[L,V]}~\Big|~{\bf L}^1([0,L],E)\right)~\geq~\sum_{i=1}^{{\bf M}_{\left(2+2^{\left(6+2/{\bf \tilde{p}}\right)}\right)\cdot {\ve\over L}}}\mathcal{M}_{2\ve}\left(\mathcal{G}^{\Psi}_{\left[L,V,2^{(4+2/{\bf \tilde{p}})}\cdot {2\ve\over L},x_i\right]}~\Big|~{\bf L}^1([0,L],E)\right)\\
~\geq~{\bf M}_{\left(2+2^{\left(6+2/{\bf \tilde{p}}\right)}\right)\cdot {\ve\over L}}\cdot \ds 2^{{\bf \tilde{p}}V\over 2\Psi\left(2^{(6+2/{\bf \tilde{p}})}\cdot {\ve\over L}\right)}~\geq~{\bf M}_{258\ve\over L}\cdot \ds 2^{{\bf \tilde{p}}V\over 2\Psi\left(256\ve\over L\right)}
\end{multline*}
and this yields the first inequality in (\ref{m1}).
\endproof
\subsection{An application to scalar conservation laws with weakly nonlinear fluxes}\label{Cons}
In this subsection, we use Theorem \ref{main} and \cite[Theorem 1]{Marc} to establish an upper bound on  the $\ve$-entropy of a set of entropy admissible weak  solutions  for a scalar conservation law in one-dimensional space
\begin{equation}\label{CL}
u_t(t,x)+f(u(t,x))_x=0\qquad\forall (t,x)\in (0,+\infty)\times \R
\end{equation}
with weakly genuinely nonlinear flux $f\in\mathcal{C}^2(\mathbb{R})$, i.e., which is not affine on any open interval such that the set 
		\begin{equation}\label{wgn}
		\{u\in\mathbb{R}~|~f''(u)\neq 0\}~~\text{is dense in}~~\mathbb{R}.
		\end{equation}
We recall that the equation (\ref{CL}) does  not
possess classical solutions since discontinuities arise in finite time even if the initial data are smooth.
Hence, it is natural to consider weak solutions in the sense of distributions
that, for the sake of uniqueness, satisfy an entropy admissibility criterion~\cite{Dafermos:Book, Kruzkov}
equivalent to the celebrated Oleinik E-condition~\cite{Oleinik}
which generalizes the classical stability conditions introduced by Lax~\cite{lax57}:

\textbf{Oleinik E-condition.} {\it 
	A shock discontinuity located at $x$ and
	connecting a left state $u^L:= u(t,x-)$ with a right state $u^R:= u(t,x+)$ is 
	entropy admissible if and only if there holds
	\begin{equation*}
	\frac{f(u^L)-f(u)}{u^L-u}~\geq~\frac{f(u^R)-f(u)}{u^R-u}
	\end{equation*}
	for every $u$ between $u^L$ and $u^R$, where  $u(t,x\pm)$ denote the one-sided limits of $u(t,\cdot)$ at $x$.}

It is well-known that the equation (\ref{CL}) generates an ${\bf L}^{1}$-contractive semigroup of solutions
$(S_t)_{t \geq 0}$ that associates, to every given initial data $u_{0} \in {\bf L}^{1}(\R) \cap {\bf L}^{\infty}(\R)$,
the unique entropy admissible weak solution $S_t u_0:= u(t,\cdot)$ of the corresponding Cauchy problem
(cfr.~\cite{Dafermos:Book, Kruzkov}). For any given $T,L,M>0$, we provide an upper bound for $\mathcal{H}_{\varepsilon}\left(S_T(\mathcal{U}_{[L,M]})\big| {\bf L}^1(\mathbb{R})\right)$ with 
\[
{\mathcal U}_{[L,M]}~:=~\Big\{ u_{0} \in {\bf L}^{\infty}(\mathbb{R}) 
\ \big| \ 
\mbox{Supp\,}(u_{0}) \subset [-L,L]\ , \ 
\| u_{0} \|_{{\bf L}^{\infty}\left(\mathbb{R}\right)} \leq M
\Big\},
\]
the set of bounded, compactly supported initial data. 

By the monotonicity of the solution operator $S_t$ and recalling
that $S_t u_0$ can be obtained as a limit of  piecewise constant front tracking
approximations~\cite[Chapter 6]{Bressan:Book}, one can show that 
\v
\begin{lemma}\label{L0} For every $L,M,T>0$ and $u_0\in  {\mathcal U}_{[L,M]}$, it holds 
	\begin{equation*}
	\big\|S_T u_0\big\|_{{\bf L}^{\infty}(\mathbb{R})}~\leq~M\qquad\mathrm{and}\qquad \mathrm{Supp}(S_T u_0)~\subseteq~\big[\!- \ell_{[L,M,T]},\, \ell_{[L,M,T]}\big]
	\end{equation*}
	where 
\begin{equation*}
		\ell_{[L,M,T]}~:=~L+T\cdot f'_M \qquad\mathrm{and}\qquad f'_M~:= ~ \sup_{|v|\leq M}~|f'(v)|\,.
		\end{equation*}
\end{lemma}
{\bf Proof.} For the proof see  \cite[Lemma 2.2]{AON4}.
\qed
\medskip

Let us introduce the function $\mathfrak d:[0,+\infty)\to [0,+\infty)$ such that 
\[
\mathfrak d(h)~=~\min_{a\in [-M,M-h]}\left(\inf_{g\in\mathcal{A}_{[a,a+h]}}\|f-g\|_{{\bf L}^{\infty}([a,a+h])}\right)
\]
with $\mathcal{A}_{[a,a+h]}$ being the set of affine functions defined on $[a,a+h]$. The  convex envelop $\Phi$ of $\mathfrak d$ is defined by 
\[
\Phi~=~\sup_{\varphi\in\mathcal{G}}\varphi\quad\mathrm{with}\quad \mathcal{G}~:=~\{\varphi:[0,+\infty)\to [0,+\infty)~|~\varphi~\mathrm{is~convex},~~\varphi(0)=0,~~\varphi\leq \mathfrak d\}.
\]
The following function 
$$
\Psi(x):=\Phi(x/2)\cdot x\qquad\forall x\in [0+\infty)
$$
is convex and satisfies the condition (\ref{Psi-cond}). As a consequence of \cite[Theorem 1]{Marc}, the following holds: 
\v
\begin{lemma}\label{L1} For any $u_0\in\mathcal{U}_{[L,M]}$, the function $S_Tu_0$ has bounded $\Psi$-total variation on $\mathbb{R}$ and 
	\begin{equation*}
	TV^{\Psi}{(S_Tu_0,\mathbb{R})}~\leq~\gamma_{[L,M,T]}:=\gamma_{[L,M]}\left(1+{1\over T}\right)
	\end{equation*}
	where $\gamma_{[L,M]}$ is a constant depending only on $L,M$ and $f$.
\end{lemma}
Recalling Corollary \ref{Eucli-cases}  for $d=1$ that for every $0<\ve \ds \leq 2L\Psi^{-1}\left({V\over 4}\right)$
	\begin{equation}\label{Eucli-cor}
	 \mathcal{H}_{\ve}\left(\mathcal{F}^{\Psi}_{[L,M,V]}~\Big|~{\bf L}^1([0,L],\mathbb{R})\right)~\leq~ \left[3\log_25+\log_2(5e)\right]\cdot {2V\over \Psi\left({\ve\over 2L}\right)}+\log_2\left({8LM\over \ve}+1\right),
	\end{equation}
we prove the following:
\medskip
	
\begin{theorem}\label{M1} Assume that $f\in \mathcal{C}^2(\mathbb{R})$ satisfies (\ref{wgn}). Then, for any  constants $L,M,T>0$, the following holds
\begin{multline*}
\mathcal{H}_{\varepsilon}\left(S_T(\mathcal{U}_{[L,M]})\Big|{\bf L}^1(\mathbb{R})\right)~\leq~\log_2\left(\frac{16M(L+T\cdot f'_M)}{\ve}+1\right)\\
+2\left[3\log_25+\log_2(5e)\right]\cdot {\gamma_{[L,M]}\left(1+{1\over T}\right)\over \Psi\left({\ve\over 4L+4T\cdot f'_M}\right)}
\end{multline*}
for every $\ve>0$ sufficiently small.
\end{theorem}
{\bf Proof.} Let us define the following set
\begin{eqnarray*}
	\tilde{S}_T(\mathcal{U}_{[L,M]})&:=&\Big\{v:\left[0,2\ell_{[L,M,T]}\right]\to [-M,M]~\Big|~\exists \ u_0\in \mathcal{U}_{[L,M]}~~ \mathrm{such~that} \\
	&~&\qquad\qquad\qquad\qquad v(x)=S_T u_0\left(x-\ell_{[L,M,T]}\right)~~\forall x\in \left[0,2\ell_{[L,M,T]}\right]\Big\}.
\end{eqnarray*}
From  Lemma \ref{L0} and Lemma \ref{L1}, it holds that
\bel{ti1}
\mathcal{H}_{\varepsilon}\left(S_T(\mathcal{U}_{[L,M]})~\Big|~{\bf L}^1(\mathbb{R})\right)~=~\mathcal{H}_{\varepsilon}\left(\tilde{S}_T(\mathcal{U}_{[L,M]})~\Big|~{\bf L}^1\left(\left[0,2\ell_{[L,M,T]}\right],~\mathbb{R}\right)\right)
\eeq
and 
\[
\tilde{S}_T(\mathcal{U}_{[L,M]})~\subseteq~\mathcal{F}^{\Psi}_{\left[2\ell_{[L,M,T]},M,\gamma_{[L,M,T]}\right]}~,
\]
where 
\[
\mathcal{F}^{\Psi}_{\left[2\ell_{[L,M,T]},M,\gamma_{[L,M,T]}\right]}~=~\Big\{g\in BV^{\Psi}\Big(\left[0,2\ell_{[L,M,T]}\right],[-M,M]\Big)~\big|~ TV^{\Psi}(g,[0, 2\ell_{[L,M,T]}])\leq \gamma_{[L,M,T]}\Big\}
\]
is defined as in Corollary \ref{Eucli-cases}. By \eqref{Eucli-cor} and \eqref{ti1}, we obtain
	\begin{multline*}
	\mathcal{H}_{\varepsilon}\left(S_T(\mathcal{U}_{[L,M]})~\Big|~{\bf L}^1(\mathbb{R})\right)=\mathcal{H}_{\varepsilon}\left(\tilde{S}_T(\mathcal{U}_{[L,M]})~\Big|~{\bf L}^1\left(\left[0,2\ell_{[L,M,T]}\right],\mathbb{R}\right)\right)\\
	~\leq~\mathcal{H}_{\varepsilon}\left(\mathcal{F}^{\Psi}_{\left[2\ell_{[L,M,T]},M,\gamma_{[L,M,T]}\right]}~\Big|~{\bf L}^1\left(\left[0,2\ell_{[L,M,T]}\right],\mathbb{R}\right)\right)\\
	~\leq~ \left[3\log_25+\log_2(5e)\right]\cdot {2\gamma_{[L,M,T]}\over \Psi\left({\ve\over 4 \ell_{[L,M,T]}}\right)}+\log_2\left({16 M \ell_{[L,M,T]}\over \ve}+1\right).
	\end{multline*}
	This completes the proof.
\endproof
\medskip

\begin{remark} 
		In general, the upper estimate of $\mathcal{H}_{\varepsilon}\left(S_T(\mathcal{U}_{[L,M]})~\Big|~{\bf L}^1(\mathbb{R})\right)$ in Theorem \ref{M1} is not optimal. 
	\end{remark}
	
We complete this subsection by  considering (\ref{CL}) with a smooth flux $f$ having polynomial degeneracy, i.e., the set $I_f=\{u\in\mathbb{R}~|~f''(u)= 0\}$ is finite and for each $w\in I_f$, there exists a natural number  $p\geq 2$ such that 
		\[
		f^{(j)}(w)~=~0\qquad\forall j\in \overline{2, p}\qquad\mathrm{and}\qquad f^{(p+1)}(w)~\neq~0.
		\]
For every $w\in I_f$, let $p_w$ be the minimal $p\geq 2$ such that $f^{(p+1)}(w)~\neq~0$. The polynomial degeneracy of $f$ is defined by 
		\[
		p_f~:=~\max_{w\in I_f} p_w~. 
		\]
Recalling \cite[Theorem 3]{Marc}, we have that $S_{T}u_0\in~BV^{1\over p_f}(\R,\R)$ and 
\[
TV^{{1\over p_f}}(S_Tu_0,\mathbb{R})~\leq~{\tilde{\gamma}_{[L,M]}} \left(1+{1\over T}\right)~=~{\tilde{\gamma}_{[L,M,T]}}
\]
for a constant $\tilde{\gamma}_{[L,M]}$ depending only on $L,M$ and $f$. This yields  
\[
\tilde{S}_T(\mathcal{U}_{[L,M]})~\subseteq~\mathcal{F}^{p_f}_{\left[2\ell_{[L,M,T]},M,\tilde{\gamma}_{[L,M,T]}\right]}~,
\]
where the set 
$$\mathcal{F}^{p_f}_{\left[2\ell_{[L,M,T]},M,\tilde{\gamma}_{[L,M,T]}\right]}~=~\left\{g\in BV^{{1\over p_f}}\left(\left[0,2\ell_{[L,M,T]}\right],[-M,M]\right)~\Big|~ TV^{{1\over p_f}}(g,[0,L])\leq \tilde{\gamma}_{[L,M,T]} \right\}$$ 
is defined as in (\ref{F-gamma}). Using \eqref{Eucli-cor} one directly obtains an extended result on the upper estimate of the $\ve$-entropy of solutions in \cite[Theorem 1.5]{AON4} for general fluxes having polynomial degeneracy.
\medskip

\begin{proposition} Assume that $f$ is smooth, having polynomial degeneracy $p_f$. Then, given the constants $L,M,T>0$, for every $\ve>0$ sufficiently small, it holds that 
		\[
		\mathcal{H}_{\varepsilon}\left(S_T(\mathcal{U}_{[L,M]})~\Big|~{\bf L}^1(\mathbb{R})\right)~\leq~{\Gamma_{[T,L,M,f]}\over \ve^{p_f}}+\log_2\left({16 (L+Tf'_M)M\over \ve}+1\right),
		\]
		where
		\[
	\Gamma_{[T,L,M,f]}~=~2^{2p_f+1}\left[3\log_25+\log_2(5e)\right]\tilde{\gamma}_{[L,M]}\left(L+T\cdot f'_M\right)^{p_f} \left(1+{1\over T}\right).
		\]
	
\end{proposition}
\begin{remark}
The above estimate is sharp in this special case. Indeed, we may exactly follow the same argument as in the proof of \cite[Theorem 1.5]{AON4} to show that 
\[
\mathcal{H}_{\varepsilon}\left(S_T(\mathcal{U}_{[L,M]})~\Big|~{\bf L}^1(\mathbb{R})\right)~\geq~\Lambda_{T,L,M,f}\cdot {1\over \ve^{p_f}}~,
\]
where $\Lambda_{T,L,M,f}>0$ is a constant depending on $L,M,T$ and $f$.
Hence, $\mathcal{H}_{\varepsilon}\left(S_T(\mathcal{U}_{[L,M]})~\Big|~{\bf L}^1(\mathbb{R})\right)$ is of the order ${1\over \ve^{p_f}}$.
\end{remark}

{\bf Acknowledgments.} This research by K. T. Nguyen was partially supported by a grant from
the Simons Foundation/SFARI (521811, NTK). The authors would like to warmly thank the anonymous referees for carefully reading the manuscript
and for their suggestions, which greatly helped in improving the paper overall.

\end{document}